\documentclass[10pt]{amsart}
\usepackage{amssymb,latexsym,epsfig}
\usepackage[numbers,sort]{natbib}

\usepackage{pstricks}
\usepackage{pst-plot}
\usepackage{pst-all}
\usepackage{pstricks-add}

  {\noindent{\bf\ref{#2}#1.}}{\vspace{\baselineskip}}

\makeatletter
\def\ulabel#1#2{\@bsphack\if@filesw {\let\thepage\relax
\def\protect{\noexpand\noexpand\noexpand}%
\xdef\@gtempa{\write\@auxout{\string
\newlabel{#1}{{#2 \@currentlabel}{\thepage}}}}}\@gtempa
\if@nobreak \ifvmode\nobreak\fi\fi\fi\@esphack} \makeatother

\newtheorem{thm}{Theorem}[section]
\newtheorem{lem}[thm]{Lemma}
\newtheorem{cor}[thm]{Corollary}

\newtheorem{question}[thm]{Question}

\theoremstyle{definition}

\newcommand{\weakly}{}
\newcommand{\nbhd}{neighborhood~}

\newcommand{\rrr}{$\mathbb{R}^3 $ }
\newcommand{\hhdf}{homotopically Hausdorff~}
\newcommand{\hhdfns}{homotopically Hausdorff}
\newcommand{\inv}{^{-1}}

\title {On small homotopies of loops}
\author{G.\ Conner}
\address{Math Department\\
         Brigham Young University\\
         Provo, UT. 84602\\
         USA}
\email{conner@math.byu.edu}
\author{M.\ Meilstrup}
\email{markhm@byu.edu}

\author{D.\ Repov\v s}
\address{Institute for Mathematics, Physics and Mechanics\\
 University of Ljubljana\\ 
 Jadranska 19, Ljubljana, Slovenia 1001. }
\email{dusan.repovs@guest.arnes.si}

\author{A.\ Zastrow}
\address{Institute of Mathematics, University of Gda\'nsk, ul.\ Wita
Stwosza 57\\ Gda\'nsk, Poland 80-952}
\email{zastrow@math.univ.gda.pl}

\author{M.\ \v Zeljko}
\address{Institute for Mathematics, Physics and Mechanics,
University of Ljubljana\\ Jadranska 19, Ljubljana, Slovenia 1001.}
\email{matjaz.zeljko@guest.arnes.si}
\thanks{This research was supported by grants
SLO-POL 12(2004--2005), BI-US/04-05/35
and P1-0292-0101-04. The first author would like to thank the Fulbright Foundation, the Ad Futura agency and the University of Ljubljana for their support.  We thank the referee for several useful suggestions and comments.}
\begin{document}

\begin{abstract}

Two natural questions are answered in the negative:\\
 \begin{itemize}
 \item ``If a space has the property that small nulhomotopic loops bound small nulhomotopies, then are loops which are limits of nulhomotopic loops themselves nulhomotopic?''
 
\item ``Can adding arcs to a space cause an essential curve to become nulhomotopic?'' \\
\end{itemize}
The answer to the first question clarifies the relationship between the notions of a space being \emph{homotopically Hausdorff} and \emph{$\pi_1$-shape injective}.
\end{abstract}

\maketitle

\pagenumbering{arabic} \setcounter{page}{1}
\section{Introduction}
Anomalous behavior in homotopy theory arises when an essential map is the uniform limit of inessential maps. Such behavior manifests itself in such oddities as pointed unions of contractible spaces being non-contractible, and (infinite) concatenations of nulhomotopic loops being essential \cite{CC2}.  
\begin{figure}[htbp]
\begin{center}

\epsfig{file=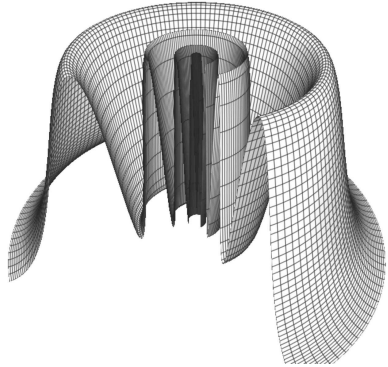,height=2in}
\caption{{ The ``surface'' portion of example $A$}}
\label{default}
\end{center}
\ulabel{innerwaved}{Figure}
\end{figure}

\begin{figure}
\begin{center}
\psset{unit=11mm,linewidth=.75pt,algebraic=true}
\begin{pspicture}(-0.5,-1.75)(12,1.75)
\psplot[plotpoints=1000]{.72}{5.5}{sin(35/(x))}
%

\psline{-}(0.5,-1.0)(0.5,1.0)%
\psline{-}(0.5,0)(3.7,0)%
\psline{-}(0.5,.5)(2.3,0.5)%
\psline{-}(0.5,-.5)(2.16,-.5)%
\psline{-}(0.5,.25)(1.6,0.25)%
\psline{-}(0.5,-.25)(1.58,-.25)%
\psline{-}(0.5,.75)(1.65,0.75)%
\psline{-}(0.5,-.75)(1.54,-.75)%

\pscircle[linewidth=1pt](9,0){1.5}
\pscircle[linewidth=1pt](9,0){1.0}
\pscircle[linewidth=1pt](9,0){0.6}
\pscircle[linewidth=1pt](9,0){0.3}
\pscircle[linewidth=1pt](9,0){0.1}
\pscircle[linewidth=1pt](9,0){0.02}

\psline{-}(7.75,0.0)(10.25,0)%
\psline{-}(9,1.25)(9,-1.25)%
\psline{-}(9.55,.55)(8.45,-.55)%
\psline{-}(8.45,.55)(9.55,-.55)%

\psline{-}(9.41,.17)(8.59,-.17)%
\psline{-}(8.59,.17)(9.41,-.17)%
\psline{-}(9.17,.41)(8.83,-.41)%
\psline{-}(9.17,-.41)(8.83,.41)%

\end{pspicture}

\caption{ Diagrams of radial projection and top view of $A$}
\end{center}
\end{figure}
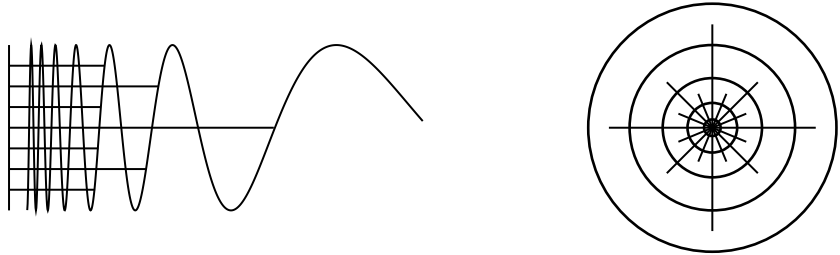

Often topologists attempt to control such behavior by requiring ``small'' maps to be nulhomotopic.  This is the flavor of the $\operatorname{k-ULC}$ property from geometric topology.

In the current article there are two natural notions of ``small'' curves which we shall study -- curves which can be homotoped into arbitarily small neighborhoods of a point, and curves which can be uniformly approximated by nulhomotopic curves.
This article describes how various embodiments of these notions are related.

In several settings one is led to ask the following

\begin{question} \ulabel{quest}{Question}
If $X$ is a space in which small nulhomotopic loops bound small homotopies, then is a loop which is the uniform limit of a family of nulhomotopic loops necessarily nulhomotopic?
\end{question}

Informally, this article is meant to clarify the above question and answer it in the negative.  Formally, this article studies two relatively new and subtly different separation axioms:  \emph{homotopically Hausdorff} and \emph{$\pi_1$-shape injective}. 

The underlying notions were introduced in a number of papers including \cite{CC2,CL,Z,CF,CC} and were put to good use in \cite{FZ1} and \cite{FZ2}.  The intuition behind a space being homotopically Hausdorff is that curves which can be homotoped into arbitrarily small neighborhoods of a point are nulhomotopic, whereas in a $\pi_1$-shape injective space one intuits that curves which can be homotoped arbitrarily close to a nulhomotopic curve are themselves nulhomotopic.
\ref{weakhhdf} shows that the property of small nulhomotopic loops bounding small nulhomotopies implies homotopically Hausdorff.

This motivates the more formal

\begin{question}
Does the homotopically Hausdorff property imply $\pi_1$-shape injectivity?
\end{question}
 
 \ref{examples} constructs two examples $A$ and $B$, neither of which is $\pi_1$-shape injective, by rotating a topologist's sine curve in \rrr to create a ``surface'' and adding a null sequence of arcs to make the space locally path connected (see the schematic diagrams for the space $A$ above).  Both spaces are homotopically Hausdorff and $B$ is \emph{strongly homotopically Hausdorff}.
 
For the sake of completeness, \ref{sinjimphh} shows shape injective implies strongly \hhdfns, and strongly \hhdf implies \weakly \hhdfns.   

The proofs that $A$ and $B$ have the desired properties require~\ref{homotopy} which answers, in the negative, the following natural
\begin{question}\ulabel{arcquest}{Question}
Can adding arcs to a space turn an essential loop into a nulhomotopic loop?
\end{question}

\begin{figure}[htbp]
\begin{center}
\ulabel{outerwaved}{Figure}
\epsfig{file=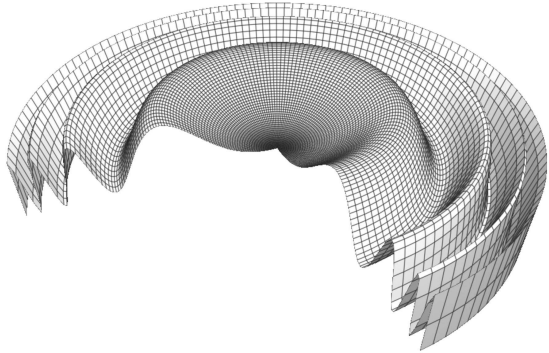,height=2in}
\caption{ The ``surface'' portion of example $B$}
\label{default}
\end{center}
\end{figure}
\subsection{Historical perspective}
In \cite{CC} the property of being \emph{homotopically Hausdorff} is described and is shown to be equivalent to the path space of the space being Hausdorff. This same notion was independently considered in \cite{Z} under the name weak $\pi_1$-continuity. The definition of \emph{shape injectivity} was introduced in \cite{CC} as the injectivity of the natural map from the fundamental group of a space into the shape group of the space.  Previously \cite{CF} studied this property extensively but it was not given a name there.  In \cite{FZ1} it is shown that shape injectivity of a space implies unique path lifting from the space to its path space and the path space is a type of \emph{generalized covering space}. 

Recently there has been renewed interest in the notion of a \emph{path space} of a separable metric space \cite{B,BS,CC,FZ1}.  The underlying desire is to find a suitable replacement for covering spaces in situations where appropriate covering spaces do not exist.  To be suitable, this replacement should have unique path lifting and so must be Hausdorff.  In \cite{CC}, the \emph{path space} is briefly described and its topology described.  In short, one uses the definition of the universal covering space in \cite{M}, but does not require the base space to be semilocally simply connected.  There are other weaker topologies which can be put on this space, see \cite{B} for instance.  However, if the path space as defined in \cite{CC} fails to be Hausdorff or has non-unique path lifting, then any weaker topology on the space will also suffer from these same deficiencies.  There is a long history of generalizations of covering spaces ranging from the work of Fox on \emph{overlays} \cite{F1,F2} up to the present \cite{MM,BS,B,CC,FZ1}.


%




\section{Definitions}
\ulabel{defs}{Section}
This section briefly recalls a number of standard definitions and introduces the definition of \emph{(strongly) homtotopically Hausdorff}.

 A \emph{curve} or $\emph{path}$ $\gamma$ in the space $X$ is a continuous function from the interval $[0,1]$ into $X$, the \emph{base} or \emph{initial point} of $\gamma$ is $\gamma(0)$, the \emph{end} or \emph{terminal point} of $\gamma$ is $\gamma(1)$,  $\gamma$ \emph{eminates} from $\gamma(0)$ and \emph{terminates} at $\gamma(1)$. Furthermore, $\gamma$ is \emph{closed} (and is called a \emph{loop}) if $\gamma(0)=\gamma(1)$, $\gamma$ is \emph{simple} if it is injective and \emph{simple closed} if it is closed and injective except at $\{0,1\}$.  Closed curves are often considered as having domain $S^1$, in the obvious way.

A \emph{free homotopy} between loops in $X$ is a continuous map from the closed annulus $[1,2] \times S^1$ to $X$ whose restriction to the boundary components of the annulus are the given loops.  A \emph{based homotopy} between loops $\gamma, \gamma'$ with the same base point is a free homotopy with the additional property that the interval $[1,2] \times\{0\}$ maps to the base point.  Two loops are \emph{(freely) homotopic} if there is a (free) homotopy between them.  A loop is \emph{nulhomotopic} or \emph{inessential} if it is homotopic to a constant map and is \emph{essential} otherwise.  A loop is nulhomotopic if and only if, when considered as a map from $S^1$ into $X$, it can be completed to a map from $B^2$ into $X$.

In \cite{CC} the \emph{path space of the space $X$ based at $x_0$}, $\Omega(X,x_0)$, is defined to be the
space of homotopy classes $\operatorname{rel}\{0,1\}$ of paths in the space $X$ based at the point $x_0 \in X$.  The path space is given the following topology: if $p$ is a path in $X$ emanating from $x_0$, and $U$ is an open neighborhood of $p(1)$, define $O(p,U)$ to be the collection of homotopy classes of paths rel$(0,1)$ containing representatives of the form $p \cdot \alpha$ where $\alpha$ is a path in $U$ emanating from $p(1),$ and take $\{O(p,U)\}$ as a basis for the topology of $\Omega(X,x_0)$.   If $X$ is semilocally simply connected then $\Omega(X,x_0)$ is the universal covering space of $X$ \cite{M}.
See \cite{CC2,FZ2} for discussions of the path space of the Hawaiian earring.

A space $X$ is \emph{$\pi_1$-shape injective} (or just \emph{shape injective}) if there is an absolute retract
$R$ which contains $X$ as a closed subspace so that whenever  $\gamma$ is an essential closed curve in $X$ then
there is a \nbhd $V$ of $X$ in $R$ such that $\gamma$ essential in $V.$ If the above
condition holds for $X$ as a closed subspace of the absolute retract
$R$, and $X$ is a closed subspace of the absolute retract $S$, then
$X$ also satisfies the above condition for $S.$ For the purposes of
this paper, $R$ will always be $\mathbb R^3.$  

If $X$ is connected, locally path connected and compact, the above definition is equivalent to the following: $X$ is $\pi_1$-shape injective if given any essential loop $\gamma$ in $\pi_1(X)$ there is a finite cover $\mathcal{U}$ of $X$ so the natural image of $\gamma$ in $\pi_1(\mathcal{N(U))}$ is essential, where $\mathcal{N(U)}$ denotes the \emph{nerve} of $\mathcal{U}.$  This is, furthermore, equivalent to the property of the natural map from $\pi_1(X)$ to \emph{shape group} of $X$, $$\displaystyle \underset{{\text{$\mathcal{U}$ a finite open cover of $X$}}}{\underleftarrow{\lim} \;\pi_1(\mathcal{N}(\mathcal{U})),}$$ is an injective homomorphism.
Thus, the name \emph{shape injective} is somewhat natural.

A space $X$ is \emph{\weakly \hhdf at a point $x_0\in X$} if for all
essential loops $\gamma$ based at $x_0$ there exists a \nbhd $U$ of $x_0$
such that no loop in $U$ is homotopic (in $X$) to $\gamma$ rel
$x_0.$ Furthermore, $X$ is \emph{\weakly \hhdf}if $X$ is \weakly \hhdf
at every point.

A space $X$ is \emph{strongly \hhdf at $x_0\in X$} if for each
essential closed curve $\gamma\in X$ there is a \nbhd of $x_0$ which
contains no closed curve freely homotopic (in $X$) to $\gamma.$ We
say that a space $X$ is \emph{strongly \hhdf}if $X$ is
strongly \hhdf at each of its points.

Intuitively, a space is \hhdf if loops which can
be made (homotopically) arbitrarily small are in fact nulhomotopic,
where the modifier \emph{strongly} allows the homotopies involved to be free homotopies.  The article \cite{CC} mentions that the name \emph{homotopically Hausdorff} was motivated by the fact that $\Omega(X,x_0)$ is Hausdorff
if and only if $X$ is homotopically Hausdorff.

Care is needed when defining these properties for non-compact spaces since we wish the notions to be topological invariants. For instance, a punctured plane is strongly \hhdf (being strongly \hhdf at each of its points), but contains an essential loop which can be homotoped to be arbitrarily small.  On the other hand, the punctured plane is homeomorphic to $S^1 \times \mathbb{R}$, endowed with its natural metric, in which no essential loop has small diameter.  

\begin{lem} \ulabel{sinjimphh}{Lemma}
\noindent \begin{enumerate}
\item If $X$ is shape injective, then $X$ is strongly \hhdfns.
\item If $X$ is strongly \hhdf at $x_0\in X$ then $X$ is \weakly \hhdf at $x_0.$
\end{enumerate}
\end{lem}

\begin{proof}
\noindent 
For part (1), suppose $X$ is a closed subspace of the absolute retract $R.$ Let
$\gamma$ be a loop in $X$ which can be freely homotoped in $X$ into an arbitrary neighborhood of $x_0$ in $X$. 
Then $\gamma$ is nulhomotopic in any \nbhd of $X$ in $R$, and since $X$ is shape
injective, $\gamma$ must be nulhomotopic. Therefore $X$ is strongly \hhdf at each of its points.

Part (2) follows immediately, since a loop which is nulhomotopic rel its base point is freely nulhomotopic.

\end{proof}

The reverse implications do not hold, not even if the space is required to be a Peano continuum. This article constructs two Peano
continua which are subspaces of $\mathbb{R}^3\!,$ $A$ and $B$,  and shows neither is shape injective while
both are \hhdf and one is even strongly \hhdfns. Both spaces will be
formed by rotating a topologist's sine curve and adding a null
sequence of arcs to make the space locally path connected.

\section{Examples}
\ulabel{examples}{Section}
The first example,  $A$, is obtained by taking the ``surface'' obtained by rotating the topologist's sine curve
about its limiting arc --a space which is not locally connected at its central arc-- and then adding a null sequence of arcs on a
countable dense set of radial cross sections to make the space
locally path connected at the central arc.  See Figures $1$ and $2.$ 

The left half
of Figure $2$ above shows a radial projection of the space $A$,
where the horizontal lines are the connecting arcs which have been added to 
the various radial cross sections. The right half of the diagram
shows a top view of the space, with the concentric circles denoting
the crests of the rotated sine curve, and the line segments
depicting the added arcs. We will refer to the various pieces of $A$
as the \emph{surface} (the rotated $\sin(1/x)$ curve), the
\emph{central limit arc}, and the \emph{connecting arcs}. Let
$\Gamma$ be the union of the interiors of the connecting arcs.

\begin{lem}\ulabel{loopb}{Lemma}
A loop $b$ of constant radius in the surface of $A$ is not freely
nulhomotopic unless it is nulhomotopic in its image.
\end{lem}

\begin{proof} If $b$ were nulhomotopic, then by \ref{homotopy} there
is a nulhomotopy of $b$ whose image does not intersect the interior
of any of the connecting arcs. Thus the image of this homotopy lies
in the path component of $b$ in $A-\Gamma$ (the complement of the
connecting arcs). This path component is the surface of $A$, which
is a punctured disc, but $b$ is not nulhomotopic in the punctured
disc unless it is nulhomotopic in its image.
\end{proof}

\begin{cor}
\ulabel{cor}{Corollary}
The Peano continuum $A$ is not shape injective.
\end{cor}

\begin{proof}
Let $b$ be a simple closed curve of constant radius in the surface, and let $a$ be a path from the base point to the initial point of $b.$ Let $\gamma = a b a\inv.$ Then $\gamma$ is nulhomotopic if and only if $b$ is, since
$b=a\inv \gamma a.$ Then, by \ref{loopb}, $\gamma$ is not
nulhomotopic, since it is not even freely nulhomotopic.

Every \nbhd of $A$ in $\mathbb{R}^3$ contains a neighborhood of the
surface union the central arc, which is just a 3-ball. Since
$\gamma$ is conjugate to a loop in the surface, $\gamma$ is then
nulhomotopic in any \nbhd of $A$, and thus $A$ is not shape
injective.
\end{proof}

\begin{thm}\ulabel{Ahhdf}{Theorem}
The space $A$ is \weakly \hhdfns, but not strongly \hhdfns.
\end{thm}

\begin{proof}
Clearly $A$ cannot be strongly \hhdf because of the loop $b$,
mentioned above, which is not nulhomotopic but can be freely
homotoped, in the surface, into any neighborhood of any point on the central limit
arc.

Now, $A$ is \weakly \hhdf at every point not in the
central arc since $A$ is locally contractible at any such point. 
In the following section, sufficient conditions for being \weakly \hhdf are proven in \ref{weakhhdf}.
Thus it remains to be shown that $A$ satisfies the hypothesis of \ref{weakhhdf} for
a base point $x_0$ in the central arc. Let $\varepsilon>0$ be given.
Choose $0<\delta<\varepsilon$, and let $\gamma$ be a nulhomotopic
loop contained in the ball $B(x_0,\delta).$ Since there are only
countably many connecting arcs, one can find a closed ball $C$ of
radius slightly greater than $\delta$ so  $\partial C$ does not intersect any
of the end points of the connecting arcs, and so  $C\subset
B(x_0,\varepsilon).$

Let $h\colon B^2\to A$ be a nulhomotopy of $\gamma.$ One can alter $h$ to obtain a nulhomotopy whose image stays in
the ball of radius $\varepsilon$ centered at $x_0.$ To do this, consider the places where the image of $h$ intersects
$\partial C$.  The following describes how to modify the homotopy so  its image remains
in a small neighborhood of $C$, which is contained in
$B(x_0,\varepsilon).$ 

First, by \ref{homotopy}, one can assume the image of $h$ does
not intersect any of the connecting arcs which $\gamma$ does not
intersect. In particular, the intersection of the image of $h$ with
$\partial C$ does not intersect any connecting arc.

Let $\ell$ denote the intersection of $C$ with the central limit
arc. Consider $h\inv(\ell)$ in $B^2.$ Let $K$ be the closure of the
union of all components of $B^2-h\inv(\ell)$ which intersect
$\partial B^2$, and let $O$ be the open set $K^c.$ 
Since $h(\partial O)\subset \ell$, and since $\ell$ is
an absolute retract, one may adjust $h$ on $O$ leaving $h|_{K}$ fixed
and sending $\overline O$ into $\ell.$ 


Thus the image of (the modified) $h$ will not pass through
$\partial C$ to exit $B(x_0,\varepsilon)$ along the central arc.

Since $C$ is a ball centered at a point on the central arc, the
intersection of $\partial C$ with the surface is a discrete collection of
circles $\{c_i\}$. Let $n$
be a component of $h\inv(c_i).$ Since 
$n$ is a component of the closed set $h\inv(c_i)$, it is closed. By
the way $C$ was chosen, $\gamma=h(\partial B^2)$ does not
intersect $\partial C$, so there exists a simple closed curve $s$
which separates $n$ from $\partial B^2.$ One can choose $s$ to be
close enough to $n$ so $s\cap h\inv(c_j)=\emptyset$ for all
$j\neq i$, and also so  $h(s)$ is contained in
$B(x_0,\varepsilon).$ Since $s$ is a simple closed curve in the disc
$B^2$, the Schoenflies theorem says $s$ bounds a disc $D$,
and then $h|_D$ is a nulhomotopy for $h(s).$ Because of the way $s$ was chosen, $h(s)$ is a loop in the surface, contained in a small
neighborhood of the circle $c_i$ which is an annulus $a_i$ contained
in $B(x_0,\varepsilon).$ Now, $c_i$ is a deformation retract of
$a_i$, and every nonzero multiple of $c_i$ is essential in $A$ by
\ref{loopb}.  If $h(s)$ were essential in $a_i$ then it would be homotopic in $A$ to a non-zero multiple of $c_i$ and would thus be essential in $A$.  Thus $h(s)$  is nulhomotopic in $a_i.$ Choose a nulhomotopy
$g_i$ for $h(s)$ which lies in $a_i$ and whose image has diameter no larger than that of $h(D)$; if $h|_D$ already lies in $a_i$, then let $g_i=h|_D$. Adjust the homotopy $h$ on the interior of the disc $D$ to be
$g_i.$ Repeat this for all components $n$ of the various preimages
$h\inv(c_i)$, to ensure the image of $h$ does not intersect
$\partial B(x_0,\varepsilon).$

The modified homotopy is still continuous, since whenever
the function was altered on a subset of $B^2$, the new image had
diameter less than or equal to the original diameter. Thus the
modified homotopy $h$ has image contained in $B(x_0,\varepsilon)$,
and so \ref{weakhhdf} applies. Therefore the space $A$ is \weakly
\hhdf at any base point $x_0$ contained in the central arc.

\end{proof}

The second example is similar to the space $A$ but is constructed by
rotating the topologist's sine curve about its initial point
$(r_0,\sin(1/r_0))$ along a vertical axis, instead of rotating about the
limit arc as depicted in Figure 3. To be precise, one can express the space $B$ in cylindrical coordinates in terms of the space $A$:
$B=\{(r_0-r,\phi,z) ~|~ (r,\phi,z)\in A, ~0\leq r\leq r_0\}.$ The central limiting arc of $A$ corresponds to the limiting \emph{outer
annulus} in $B$, and the \emph{connecting arcs} of $B$ limit on every
point of this annulus. The \emph{surface} of $B$ is
homeomorphic to $\mathbb R^2.$ Again, let $\Gamma$ denote the union
of the interiors of the connecting arcs.

\begin{thm}
The Peano continuum $B$ is not shape injective.
\end{thm}

\begin{proof}
Consider a loop about the limiting outer annulus. As in the proof of
\ref{loopb}, one applies \ref{homotopy} to see that this loop
cannot be nulhomotopic, as it is not nulhomotopic in the path component
of $B-\Gamma$ containing it, which is an annulus. Any \nbhd of $B$
in $\mathbb R^2$ contains a \nbhd of the annulus which intersects
the rotated surface, and this curve can then be homotoped into the
surface. Thus this curve is nulhomotopic in any \nbhd of $B.$ Once again,
one conjugates by a path to the base point. Since the original
path is essential it follows that $B$ is not shape injective.

\end{proof}

\begin{thm}
\ulabel{bthm}{Theorem}
The space $B$ is strongly \hhdfns.
\end{thm}

\begin{proof}
This proof will be remarkably similar to the proof of \ref{Ahhdf}.
Evidently $B$ is strongly \hhdf at any point outside the outer
annulus, since it is locally contractible there. Let $x_0$ be a point on
the outer annulus.

Let $\varepsilon>0$ and choose $0<\delta<\varepsilon.$ Let
$\gamma,\gamma'$ be homotopic essential loops contained in the ball
$B(x_0,\delta)$. One can find a
closed ball $C$ contained in $B(x_0,\varepsilon)$ and containing $B(x_0,\delta)$ so  $\partial C$ does
not intersect any of the end points of the connecting arcs.

Let $h\colon (S^1\times I)\to B$ be a homotopy between $\gamma$ and
$\gamma'.$ By \ref{homotopy}, one may assume  the image of $h$
does not intersect any of the connecting arcs which do not intersect
either $\gamma$ or $\gamma'.$ In particular, the intersection of the image of $h$ with
$\partial C$ does not intersect any connecting arc.

Let $d$ denote the disc which is the intersection of $C$
with the outer annulus. Similar to the proof of \ref{Ahhdf}, since
$d$ is an absolute retract ($d$ plays the role of $\ell$ in that
proof), $h$ may be altered so every component of $(S^1\times
[0,1])-h\inv(d)$ intersects the boundary of $S^1\times[0,1]$. Thus the image of (the
modified) $h$ cannot pass through $\partial C$ to exit
$B(x_0,\varepsilon)$ along the outer annulus.

Since $C$ is a ball centered on the outer annulus, the intersection
of $\partial C$ with the surface is a discrete collection of circles
$\{c_i\}$, each of which bounds a disc $d_i$ in the surface. Let $n$ be a component of
$h\inv(c_i).$ Since $h\inv(c_i)$ is closed, $n$ is
closed. Because of the way $C$ was chosen, $\gamma\cup\gamma'
=h(S^1\times\{0,1\})$ does not intersect $\partial C.$ If
$h\inv(c_i)$ separates $S^1\times\{0\}$ from $S^1\times\{1\}$, then
both $\gamma,\gamma'$ will be nulhomotopic, since they are both
homotopic to a power of $c_i$ which bounds the disc $d_i$,
contradicting the choice of an essential curve $\gamma.$

Since $h\inv(c_i)$ does not separate $S^1\times\{0,1\}$, there
exists a simple closed curve $s$ which separates $n$ from
$S^1\times\{0,1\}.$ Choose $s$ to be close enough to $n$ so $s\cap h\inv(c_j)=\emptyset$ for all $j\neq i$, and also so
 $h(s)$ is contained in $B(x_0,\varepsilon).$

By the way $s$ was chosen, $h(s)$ is a loop in the surface, which lies
in a small \nbhd of $c_i$, and thus bounds a disc in a small \nbhd
of the disc $d_i.$ Let $g$ denote a nulhomotopy of $h(s)$ whose image lies in this disk and has diameter no larger than that of $h(s).$  By the
Schoenflies theorem, the component of the complement of $s$ which
does not contain the boundary $S^1\times\{0,1\}$ is a disc. Adjust
the homotopy $h$ on this disc to be the nulhomotopy $g$; $h$ need not be adjusted in the case that it's image is already sufficiently small.  Carry out this adjustment for all components $n$ of the various preimages $h\inv(c_i)$,
to ensure  the image of $h$ does not intersect $\partial
B(x_0,\varepsilon).$

As in the proof of \ref{Ahhdf}, the modified homotopy is still
continuous, since whenever the function was altered on a subset of
$S^1\times[0,1]$ the new image had diameter less than
or equal to the original diameter. Thus the modified homotopy $h$
has image contained in $B(x_0,\varepsilon)$, and so \ref{stronghhdf}
applies. Therefore the space $B$ is strongly \hhdf at any point
$x_0$ contained in the outer annulus.

\end{proof}

\section{Technical Lemmas}
\ulabel{lemmas}{Section}
First, a condition on metric spaces is given which implies the condition \hhdf at a
point.  The basic idea
is that for every small nulhomotopic curve, there is a nulhomotopy
of small diameter. This is similar to $\operatorname{1-ULC}$, but the
condition is only required to hold for nulhomotopic loops.

\begin{lem}\ulabel{weakhhdf}{Lemma}
Suppose the metric space $X$ contains a point $x_0$ enjoying the
property that for every $\varepsilon>0$ there is a $\delta>0$ such
that for every continuous function $f\colon B^2\to X$ with $f(S^1)\subset
B(x_0,\delta)$, there is a continuous function $g\colon B^2\to X$ such that
$g|_{S^1}=f|_{S^1}$, and
$g(B^2)\subset B(x_0,\varepsilon).$ %
Then $X$ is \weakly \hhdf at $x_0.$
\end{lem}

\begin{proof}
Let $\gamma_i$ be a null sequence of loops based at $x_0$
representing the same homotopy class in $\pi_1(X,x_0).$ Construct
a nulhomotopy $f$ of $\gamma_1$ as follows: Consider a Hawaiian
earring in the disc $B^2$ as in Figure 4 below. Define $f$ on each
of the arcs $c_i$ to be the loops $\gamma_i$ in $X.$ Then each
portion $D_i$ of the disc where $f$ is not yet defined is bounded by
a curve $\gamma_i \overline{\gamma}_{i+1}$, which is nulhomotopic.
Thus $f$ can be defined on the entire disc so that it is
continuous at every point except possibly at the base point of the
Hawaiian earring.



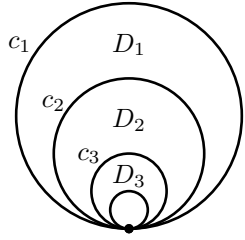
\begin{figure}
\begin{center}
\psset{unit=1cm,linewidth=.75pt,algebraic=true}
\begin{pspicture}(-2,-.5)(2,3.5)

\pscircle[linewidth=1pt](0,1.5){1.5}
\pscircle[linewidth=1pt](0,1.0){1.0}
\pscircle[linewidth=1pt](0,0.5){0.5}
\pscircle[linewidth=1pt](0,0.25){0.25}

\pscurve[showpoints=true](0,0)(0,0)

\rput(-1.45,2.45){\Rnode{e1}{$c_1$}}
\rput(-1.0,1.67){\Rnode{e2}{$c_2$}}
\rput(-0.53,0.97){\Rnode{e3}{$c_3$}}

\rput(0.0,2.45){\Rnode{d1}{$D_1$}}
\rput(0.0,1.45){\Rnode{d2}{$D_2$}} \rput(0.0,0.7){\Rnode{d3}{$D_3$}}

\end{pspicture}
\caption{The Hawaiian earring}
\end{center}
\end{figure}

One carefully chooses nulhomotopies $f|_{D_i}$ of $\gamma_i
\overline{\gamma}_{i+1}$ to ensure the continuity of $f$
at the base point. Let $(\epsilon_n)$ be a sequence of positive
numbers decreasing to 0. 
Then by hypothesis there exists a sequence $(\delta_n)$ such that any
nulhomotopic loop contained in the ball $B(x_0,\delta_n)$ has a
nulhomotopy whose image is contained in the ball $B(x_0,
\epsilon_n).$ Without loss of generality, assume
$\delta_n\geq \delta_{n+1}.$ Since the loops $\gamma_i$ form a null
sequence limiting to a point, choose $k_n$ to be the minimal
index so $\gamma_{i} \overline{\gamma}_{i+1}$ has diameter
less than $\delta_n$ for all $i\geq k_n.$ Then since $\delta_n\geq
\delta_{n+1}$, it follows that $k_n\leq k_{n+1}.$ Then for all $k_n\leq
i <k_{n+1}$, define $f|_{D_i}$ to be a nulhomotopy of $\gamma_i
\overline{\gamma}_{i+1}$ with diameter less than $\epsilon_n$, which
exists by hypothesis.

To see this defines $f$ as a continuous function at the base
point $y$ of the Hawaiian earring, let $\varepsilon>0$ be given.
Then there is some $n$ such that $\epsilon_n<\varepsilon.$ Then
by the construction, the arc $c_{k_n}$ in the disc (which maps to
$\gamma_{k_n}$) bounds a disc whose image is contained in $B(x_0,\epsilon_n)\subset B(x_0,\varepsilon).$
 Since there are only finitely many
discs $D_i$, for $i<n$, one can find a $\delta>0$ such that $f$ maps
$\left(\bigcup_{i=1}^{n-1} D_i\right) \cap B(y,\delta)$ into
$B(x_0,\varepsilon).$ Then since $f$ maps $\bigcup_{i=n}^\infty D_i$
into $B(x_0,\varepsilon)$, it follows that $f(B(y,\delta))\subset
B(x_0,\varepsilon)$, and thus $f$ is continuous.

Thus $\gamma_1$ is nulhomotopic and consequently all of the curves $\gamma_i$ are nulhomotopic.  Therefore
$X$ is \weakly \hhdf at $x_0.$
\end{proof}

While this condition is sufficient, it is not necessary.
Consider the cone over the Hawaiian earring in \ref{hecone}, which is contractible,
hence \hhdfns. The loops of the base Hawaiian earring are
nulhomotopic, yet they require nulhomotopies of large diameter
(passing over the cone point).

\begin{figure}
\begin{center}

\epsfig{file=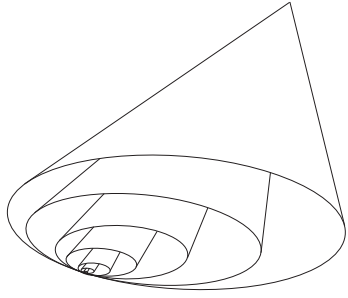,height=2in}

\caption{A space which is homotopically Hausdorff but which does not satisfy the conditions of \ref{weakhhdf}}
\ulabel{hecone}{Figure}

\label{default}
\end{center}
\end{figure}

We now describe a condition which is sufficient to imply strongly \hhdf at a
point; it guarantees a homotopy of small diameter between every pair of essential homotopic
curves nearby a point.

\begin{lem}\ulabel{stronghhdf}{Lemma}
%
Let $X$ be a metric space and $x_0\in X$ such that for every
$\varepsilon>0$ there is a $\delta>0$ such that for every essential map of an annulus
$f\colon S^1\times[0,1]\to X$ with $f(S^1\times\{0,1\})\subset B(x_0,\delta)$,
there is a map of an annulus $g\colon S^1\times[0,1] \to X$ such that
$g|_{S^1\times\{0,1\}}=f|_{S^1\times\{0,1\}}$, and
$g(S^1\times[0,1])\subset B(x_0,\varepsilon).$ %
Then $X$ is strongly \hhdf at $x_0.$
\end{lem}

\begin{proof}
Let $\gamma_i$ be a null sequence of loops which are freely
homotopic to each other and which converge to $x_0.$ It must be shown that
 $\gamma_1$ is freely nulhomotopic. Suppose not. Then each $\gamma_i$ is essential. By way of contradiction, one
constructs a nulhomotopy $f$ of $\gamma_1$ as follows. In the unit
disc $B^2$, specify concentric circles $c_i$ of radius $1/i.$ Define
$f|_{c_i}=\gamma_i$ in $X.$ Let $A_i$ be the annulus in $B^2$
bounded by $c_i \cup c_{i+1}.$ Since $f|_{c_i}=\gamma_i$, and since
$\gamma_i$ is homotopic to $\gamma_{i+1}$, one can extend $f$ to each $A_i$ in such a way that, by an argument similar to the end of \ref{weakhhdf}, we may extend $f$ to a map which is also continuous at the center point of $B^2.$ %
Thus $f$ is a nulhomotopy of the curve $\gamma_1$, and hence
$X$ is strongly \hhdf at $x_0.$
\end{proof}

The next lemma can be thought of as a general position result for arcs and nulhomotopies.  It says that if a nulhomotopic loop does not meet the interiors of a collection of arcs, then there is nulhomotopy for the loop which does not meet the interiors of the arcs.

\begin{lem}\ulabel{homotopy}{Lemma}
Let $X$ be a topological space. Let $\Xi$ be a disjoint union of open
sets in $X$, each of which is homeomorphic to an open arc,  and let
$Z=X-\Xi.$
\begin{enumerate}

\item \label{lab1}Let $g\colon B^2 \to X$ be a nulhomotopy such that
$g(\partial B^2) \subseteq Z.$   Then there is a nulhomotopy $g'\colon B^2 \to Z$
with $g\mid_{\partial B^2}=g'\mid_{\partial B^2}.$

\item \label{lab2} Let $h\colon (S^1\times[0,1])\to X$ be a homotopy between
two essential curves $\gamma$ and $\gamma'$ in $Z.$ Then there is a
homotopy $h'$ between $\gamma$ and $\gamma'$ such that the
image of $h'$ lies in $Z.$
\end{enumerate}
\end{lem}
An alternate way of stating the conclusion of this theorem would be to say the natural map $i_*:\pi_1(Z) \to \pi_1(X)$ induced by inclusion is injective.
\begin{proof}
For each arc $\xi$ in $\Xi$, let $a_\xi$ be an open arc in
$\xi$ whose closure is contained in $\xi$, and let $A$ be the
union of the arcs $a_\xi.$  The subspace $Z$ is a
strong deformation retract of $X-A$, so it suffices to show that the
maps described above take values in $X-A.$

For the moment we proceed with the proof of (\ref{lab2}). Let $K$ be
the boundary of the component of $S^1\times[0,1]-h^{-1}(A)$
containing $S^1\times\{0\}.$  Now, $h$ is constant on each
component of $K$ since the boundary of $A$ is totally disconnected.

Suppose $K$ separates $S^1\times\{0,1\}$, the boundary of the
annulus.  Since $\mathbb{R}^2$ is \emph{unicoherent}, the interior
of the annulus $S^1\times(0,1)$ has the following property: if a
compact subspace contained in the interior of the annulus separates
two points of the closed annulus, then one of the components of the
subspace separates those two points.  Consequently, one may choose a
component, $T$, of $K$ which separates $S^1\times\{0\}$ from
$S^1\times\{1\}.$  One now creates a new map which is equal to
$h$ everywhere except for the component of the complement of $T$
which contains $S^1\times\{1\}$ and defines it to be the constant
$h(T)$ on that component.  Furthermore, one may adjoin a disk to
$S^1\times[0,1]$ along $S^1\times\{1\}$ to obtain $B^2$ and extend
the new map by defining it to be the constant $h(T)$ on this disk
also.  The result is a nulhomotopy of $\gamma$, contradicting the hypothesis that $\gamma$ is essential.


Thus assume  $K$ does not separate $S^1\times\{0,1\}.$  Recall $h$ is constant on each component of $K.$ Hence we may define $h'$ by having it agree with $h$ on the component of the complement $K$ which contains $S^1\times\{0,1\}$ and defining it to be constant on the other components of the complement of $K$, thus proving (\ref{lab2}).

To prove (\ref{lab1}) it is enough to mention that if  $L$ is the boundary of the component  of $B^2-g^{-1}(A)$ containing $\partial B^2$, then, as in the argument above, $g$ is constant on each component of $L.$  Define $g'$ to be equal to $g$ on the component of the complement of $L$ containing $\partial B^2$ and to be constant on the other components of the complement of $L.$
\end{proof}

\end{document}